\documentclass[10pt]{article}
\usepackage{amssymb, amsmath, amsfonts}
\usepackage[cp1251]{inputenc}
\usepackage[english,russian]{babel}
\usepackage[mathscr]{eucal}
\usepackage{graphics}

\newtheorem{theorem}{Теорема}
\newtheorem{definition}{Определение}

\title{Вокруг теоремы А.Д. Александрова\\ о характеризации сферы}
\author{В.А. Александров\\
 \\
\small\textit{Институт математики им. С.~Л.~Соболева СО РАН},\\
\small{\textit{пр. академика Коптюга 4, 630090, Новосибирск, Россия} и}\\
\small\textit{физический факультет Новосибирского государственного университета,}\\
\small\textit{ул. Пирогова 2, 630090, Новосибирск, Россия}\\
\small{e-mail: alex@math.nsc.ru}}
\date{20 декабря 2012 года}

\begin{document}
\maketitle
\begin{abstract}
\textbf{Victor Alexandrov,} 
\textit{Around the A.D. Alexandrov's theorem on a characterization of a sphere}.
This is a survey paper on various results relates 
to the following theorem first proved by A.D. Alexandrov: 
\textit{Let $S$ be an analytic convex 
sphere-homeomorphic surface in $\mathbb R^3$ and 
let $k_1(\boldsymbol{x})\leqslant k_2(\boldsymbol{x})$ be 
its principal curvatures at the point $\boldsymbol{x}$. 
If the inequalities 
$k_1(\boldsymbol{x})\leqslant k\leqslant k_2(\boldsymbol{x})$ 
hold true with some constant $k$ for all 
$\boldsymbol{x}\in S$ then $S$ is a sphere.}
The imphases is on a result of Y. Martinez-Maure who first proved that
the above statement is not valid for convex $C^2$-surfaces.
For convenience of the reader, in addendum we give 
a Russian translation of that paper by Y. Martinez-Maure
originally published in French in 
\textit{C. R. Acad. Sci., Paris, S\'{e}r. I, Math.} 
{\bf 332} (2001), 41--44.
\end{abstract}

\renewcommand{\thefootnote}{\fnsymbol{footnote}}

\section{Введение}
\footnotetext{Работа выполнена при поддержке РФФИ (код проекта 
10--01--91000--анф),
ФЦП <<Научные и научно-педагогические кадры инновационной России>> на 
2009--2013 гг. (гос. контракт 02.740.11.0457) и Совета по грантам 
Президента РФ для поддержки ведущих научных школ (НШ--6613.2010.1).}
Эта статья является расширенной версией одноимённого
доклада, прочитанного автором в сентябре 2011 года
на конференции <<Дни геометрии в Новосибирске>>.
Она содержит краткий обзор результатов, 
связанных со следующей теоремой
А.Д. Александрова о характеризации сферы:
\begin{theorem}\label{1}
Пусть $S$ --- аналитическая выпуклая поверхность в $\mathbb{R}^3$,
гомеоморфная сфере, и пусть 
$0\leqslant k_1(\boldsymbol{x})\leqslant k_2(\boldsymbol{x})$ ---
её главные кривизны в точке $\boldsymbol{x}$. 
Если существует постоянная $k$ такая, что для всех
$\boldsymbol{x}\in S$ выполняются неравенства 
$ k_1(\boldsymbol{x})\leqslant k\leqslant k_2(\boldsymbol{x})$, то
поверхность $S$ является сферой.
\end{theorem}

В разделе 2 мы расскажем о связях теоремы \ref{1} с задачами однозначной
определённости выпуклых поверхностей, в разделе 3 --- о связях с теоремами 
о поверхностях Вейнгартена. 
В разделе 4 мы обсудим особенно важную для нашего 
изложения работу И. Мартинеса-Мора \cite{MM01}, показавшего, 
что для выпуклых $C^2$-поверхностей утверждение, аналогичное теореме \ref{1}, 
не верно. В разделе 5 мы коснёмся статьи \cite{Pa05} Г.Ю. Паниной, 
продолжившей исследования И. Мартинеса-Мора (и впервые показавшей, 
что для выпуклых $C^2$-поверхностей утверждение, аналогичное 
теореме \ref{1},  не верное), а в разделе 6
говорим о перспективах дальнейших исследований, связанных с 
тематикой настоящей статьи.
Наконец, для удобства читателя в <<Приложении>>
мы приводим перевод статьи И. Мартинеса-Мора \cite{MM01}.

\section{Теоремы единственности для выпуклой поверхности\\ 
с заданной функциональной зависимостью между\\ 
главными кривизнами и вектором нормали}

Теорема \ref{1} является очевидным частным случаем теоремы 1-а из
статьи А.Д. Александрова 1966 года \cite{AD66}. Сама эта теорема 1-а 
сформулирована в терминах так называемой относительной дифференциальной геометрии
(именно поэтому мы не приводим здесь её формулировку). 
В \cite{AD66} дано доказательство теоремы 1-а 
и показано, что она эквивалентна следующей теореме (см. \cite[теорема 1]{AD66}):

\begin{theorem}\label{AD66-1}
Пусть $S$ и $S^0$ --- аналитические выпуклые поверхности в $\mathbb{R}^3$,
гомеоморфные сфере, причём гауссова кривизна поверхности $S^0$ всюду положительна
(а значит, сама поверхность $S^0$ выпуклая).
Утверждается, что либо $S$ получена из $S^0$ параллельным 
переносом, либо существуют такие точки
$\boldsymbol{x}\in S$, $\boldsymbol{x}^0\in S^0$ 
с параллельными и сонаправленными нормалями, 
что кривизна любого нормального 
сечения в $\boldsymbol{x}$ отлична от кривизны параллельного 
нормального сечения в $\boldsymbol{x}^0$.
\end{theorem}

Убедимся, что теорема \ref{1} является следствием теоремы \ref{AD66-1}.
Пусть для поверхности $S$ выполнены условия теоремы \ref{1}.
В качестве поверхности $S^0$ возьмём сферу радиуса $1/k$.
Если $S$ не получена из $S^0$ параллельным переносом, 
то по теореме \ref{AD66-1} найдётся точка
$\boldsymbol{x}\in S$, в которой кривизна любого нормального сечения 
не равна $k$.
Но в точке $x$ имеется нормальное сечение с кривизной $k_1(\boldsymbol{x})$
и есть нормальное сечение с кривизной $k_2(\boldsymbol{x})$. Кроме того, кривизна 
нормального сечения, проходящего через точку $\boldsymbol{x}$, есть функция непрерывная, 
а $k_1(\boldsymbol{x})\leqslant k\leqslant k_2(\boldsymbol{x})$.
Значит, в точке $\boldsymbol{x}$ найдётся нормальное сечение кривизны $k$.
Полученное противоречие показывает, что поверхность $S$ 
получена из $S^0$ параллельным переносом, 
т.\,е. является сферой. Тем самым мы вывели теорему \ref{1}
из теоремы \ref{AD66-1}.

В свою очередь теорема \ref{AD66-1} была сформулирована и 
доказана А.Д. Александровым в 1966 г. в \cite{AD66} для того, 
чтобы получить новое доказательство следующей теоремы
(см. \cite[теорема 2]{AD66}):

\begin{theorem}\label{AD66-2}
Пусть $S$ и $S^0$ --- такие поверхности, как в теореме \ref{AD66-1},
и пусть $k_1(\boldsymbol{x})\leqslant k_2(\boldsymbol{x})$,  
$k_1^0(\boldsymbol{x}^0)\leqslant k_2^0(\boldsymbol{x}^0)$
их главные кривизны в точках $\boldsymbol{x}\in S$ и 
$\boldsymbol{x}^0\in S^0$.
Пусть $f(\xi, \eta, \boldsymbol{n})$ ---
такая функция численных переменных $\xi$, $\eta$ и 
единичного вектора  $\boldsymbol{n}\in\mathbb{R}^3$, что
при  $\xi>\xi'$ и $\eta>\eta'$ всегда
$f(\xi, \eta, \boldsymbol{n})>f(\xi', \eta', \boldsymbol{n})$.
Тогда если для любой точки $\boldsymbol{x}\in S$ и соответствующей ей
точки $\boldsymbol{x}^0\in S^0$, имеющей такой же вектор нормали
$($т.\,е. такой, что 
$\boldsymbol{n}^0(\boldsymbol{x}^0)=\boldsymbol{n}(\boldsymbol{x}))$
выполняется равенство
$$f(k_1(\boldsymbol{x}), k_2(\boldsymbol{x}), 
\boldsymbol{n}(\boldsymbol{x}))=
f(k_1^0(\boldsymbol{x}^0), k_2(\boldsymbol{x}^0), 
\boldsymbol{n}^0(\boldsymbol{x}^0)),$$
то поверхность $S$ получена из $S^0$ параллельным переносом.
\end{theorem}

Легко понять, что теорема  \ref{AD66-2} является простым
следствием теоремы \ref{AD66-1}.
В самом деле, если выполнены условия теоремы \ref{AD66-2},
то, ввиду монотонности функции $f$, в каждой точке
$\boldsymbol{x}\in S$ либо 
$k_1(\boldsymbol{x})\geqslant k_1^0(\boldsymbol{x}^0)$,
$k_2(\boldsymbol{x})\leqslant k_2^0(\boldsymbol{x}^0)$, либо
$k_1(\boldsymbol{x})\leqslant k_2(\boldsymbol{x}^0)$,
$k_2(\boldsymbol{x})\geqslant k_2(\boldsymbol{x}^0)$.
Поэтому в соответствующих точках $\boldsymbol{x}$, 
$\boldsymbol{x}^0$
всегда есть пара параллельных нормальных сечений с
равными кривизнами. 
А тогда по теореме \ref{AD66-1} поверхность $S$
получена из $S^0$ параллельным переносом.

Впервые утверждение, аналогичное теореме \ref{AD66-2},
было сформулировано и доказано А.Д. Александровым в 1938 г.  
\cite[теорема I]{AD38} в качестве общей теоремы единственности,
охватывающей классические теоремы Г. Минковского и Э.Б. Кристоффеля
\cite[\S 94, 95]{BF34}.
Как известно, последние утверждают, что выпуклая поверхность однозначно
определяется заданием своей гауссовой (соответственно,
средней) кривизны как функции от единичной нормали к поверхности.
Чтобы получить теорему Г. Минковского (или Э.Б. Кристоффеля)
достаточно применить теорему \ref{AD66-2}, выбрав в качестве
$f(\xi,\eta,\boldsymbol{n})$  функцию  $\xi\cdot\eta$ 
(соотв., $\xi+\eta$).

Теоремам единственности для выпуклых поверхностей 
(гладких, многогранных, многомерных, в неевклидовых 
пространствах и т.\,д.) с заданной функциональной
зависимостью между главными кривизнами и вектором единичной
нормали к поверхности посвящено большое количество работ.
При этом используются самые разные подходы и методы.
Пожалуй можно сказать, что основным методом является
свед\'{е}ние геометрической проблемы к решению некоторой
задачи для уравнения в частных производных. 
Читатель может подробнее познакомиться 
с этой проблематикой, например, по классическому
учебнику \cite{BVK73} или недавней статье \cite{GM03} 
и указанной в них литературе.

\section{Поверхности Вейнгартена}

Теорема \ref{1} связана ещё и с другим важным классом
поверхностей --- поверхностями Вейнгартена. 
Напомним, что поверхность называют поверхностью Вейнгартена,
если её главные кривизны $k_1(\boldsymbol{x})$, 
$k_2(\boldsymbol{x})$ всюду связаны 
некоторым функциональным соотношением, т.\,е. если 
$f(k_1(\boldsymbol{x}),k_2(\boldsymbol{x}))\equiv 0$, где $f$ --- некоторая 
невырожденная функция двух вещественных переменных.
При этом, в отличие от изложенного в предыдущем разделе,
не требуется выпуклости поверхности и функциональное 
соотношение, связывающее главные кривизны, не 
зависит от вектора нормали.

Эти поверхности введены Ю. Вейнгартеном 
(Julius Weingarten, 1836--1910) в связи с изучением вопроса
об описании всех поверхностей, изометричных данной.
Хотя этот вопрос в полном объёме не решён до сих пор,
поверхности Вейнгартена (наряду с поверхностями вращения)
служат важным источником примеров, на которых 
проверяются общие теории, и который привлекает внимание
специалистов  \cite[Problem 58]{Ya82}. 
Не углубляясь в детали отметим только, что частными случаями
поверхностей Вейнрагтена являются поверхности постоянной
гауссовой кривизны, минимальные поверхности и поверхности
постоянной средней кривизны.
Читатель может ближе познакомиться с 
поверхностями Вейнгартена, например, по учебнику 
\cite[Appendix to Chap. 2]{To06},
монографии \cite[Chap. V, Sect. 2]{Ho83}, или сравнительно 
недавней статье \cite{GM10} и указанной там литературе.

Поверхности Вейнгартена обладают многими замечательными 
свойствами. Одним из наиболее глубоких является 
следующая теорема Х. Хопфа (Heinz Hopf, 1894--1971;
см. \cite{Ho51} или \cite[Chap. VI, Sect. 2]{Ho83}):

\begin{theorem}\label{Ho51}
Пусть $S$ является погруженной в $\mathbb R^3$
замкнутой поверхностью рода $0$, имеющей постоянную среднюю
кривизну. Тогда $S$ является сферой.
\end{theorem}

У этой теоремы довольно длинная история, 
входить во все детали которой мы здесь не будем. 
Отметим только, что для случая выпуклых поверхностей она 
была впервые доказана Г. Либманом %Heinrich Liebmann (1874-1939)
в 1900 году \cite{Li00}, а аналоги и обобщения 
постоянно появляются и сейчас (см, например, \cite{AR04}).

Легко видеть, что теорема \ref{Ho51} является 
очевидным следствием из теоремы \ref{AD66-2}.
Для этого нужно взять в качестве $S^0$ сферу и положить
$f(k_1,k_2)=k_1+k_2$.
Аналитичности поверхности $S$ в теореме \ref{Ho51}
требовать не нужно, поскольку для поверхности
постоянной средней кривизны она имеется <<автоматически>>.

Для изучения поверхностей Вейнгартена разные авторы применяли
разные методы. Пожалуй, основными являются 
построение голоморфных квадратичных дифференциалов и
оценивание индекса
линий кривизны поверхности в изолированной омбилической точке.
Именно последним методом в 1967 году Г.Ф. Мюнцнер \cite {Mu67}
установил многие свойства поверхностей постоянной средней кривизны
и, в частности, передоказал теоремы \ref{AD66-1} и \ref{Ho51}.

Следующий (по времени) принципиальный результат, связанный
с теоремой~\ref{1}, был опубликован в 2001 году И. Мартинесом-Мором
\cite{MM01}. Для его изложения нам понадобится понятие ежа,
которое и вводится в следующем разделе.

\section{Гладкие ежи}

Понятие гладкого ежа служит обобщением понятия 
гладкой выпуклой поверхности в $\mathbb R^n$ (точнее ---
границы гладкого выпуклого тела).

В самом деле, гладкую выпуклую поверхность можно рассматривать как
огибающую гладкого семейства опорных плоскостей данного выпуклого тела
(рис. 1, $a$). При этом с аналитической точки зрения 
семейство опорных плоскостей 
$$
\{\boldsymbol{x}\in\mathbb R^n\vert 
(\boldsymbol{x},\boldsymbol{p})=h(\boldsymbol{p})\}\eqno(1)
$$
задаётся посредством гладкого отображения 
$h:\mathbb S^{n-1}\to\mathbb R$, 
сама выпуклая поверхность допускает параметризацию
$x_h:\mathbb S^{n-1}\to\mathbb R^n$ с помощью формулы
$$
\boldsymbol{p}\mapsto (\mbox{grad\,}h)(\boldsymbol{p})
+h(\boldsymbol{p})\boldsymbol{p},\eqno(2)
$$ 
а гауссова кривизна поверхности в точке $x_h(\boldsymbol{p})$
оказывается равной определителю касательного отбражения $x_h$ 
в точке $\boldsymbol{p}\in\mathbb{S}^{n-1}$.

\begin{figure}\center
 \includegraphics{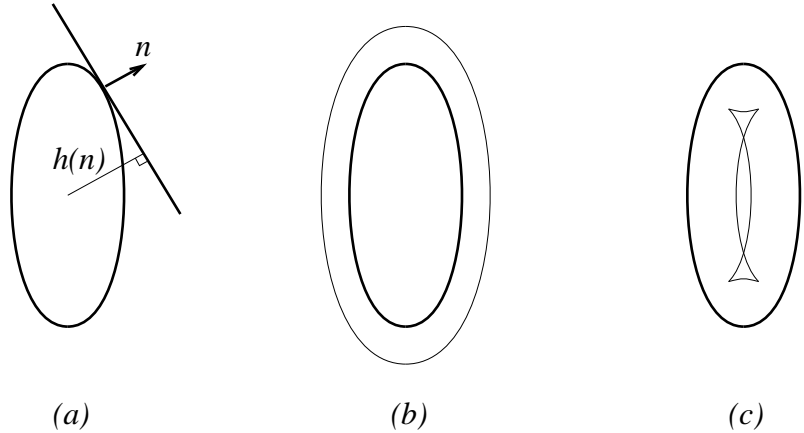}
 \caption{$(a):$ Опорная функция $h(n)$ выпуклого тела.\newline 
 $(b):$  Огибающая опорной функции $h(n)+\rm{const}$ при $\rm{const}>0$.\newline 
 $(c):$ Огибающая опорной функции $h(n)+\rm{const}$ при $\rm{const}<0$.} 
\end{figure}

Приняв такую аналитическую точку зрения можно 
рассматривать огибающие семейства плоскостей (1),
соответствующих произвольной гладкой функции
$h:\mathbb S^{n-1}\to\mathbb R$, 
независимо от того является $h$ опорной функцией
некоторого выпуклого тела или нет.
При этом огибающая может иметь особенности.
Например, из рис. 1, $b$ 
и рис. 1, $c$ видно,
что если  $h$ --- опорная функция гладкого выпуклого тела, 
то огибающая семейства плоскостей (1), соответствующих 
функции $h+\mbox{const}$, будет гладкой поверхностью,
если  $\mbox{const}>0$, но может иметь особенности,
если $\mbox{const}<0$.

Эти рассуждения приводят нас, вслед за 
Р. Ланжевеном, Г. Левиттом и Г. Розенбергом \cite{LLR88},
к следующему определению (см. также \cite{MM01, MM10})

\begin{definition}
Гладким ежом в $\mathbb{R}^n$, соответствующим
гладкой функции   $h:\mathbb S^{n-1}\to\mathbb R$, 
называется огибающая семейства плоскостей $(1)$,
т.\,е. образ отображения $x_h:\mathbb S^{n-1}\to\mathbb R^n$,
заданного формулой $(2)$.
При этом исходная функция $h$ называется опорной функцией ежа,
собственные значения касательного отображения 
$T_{\boldsymbol{p}}x_h$ в точке 
$\boldsymbol{p}\in\mathbb{S}^{n-1}$
называются главными радиусами кривизны ежа в точке
$x_h(\boldsymbol{p})$,
их произведение $($т.\,е. $\det(T_{\boldsymbol p}x_h))$
называется функцией кривизны $R_h$ ежа,
а величина $1/R_h$ --- гауссовой кривизной ежа.
Ежа, соответствующего
гладкой функции $h:\mathbb S^{n-1}\to\mathbb R$,
обозначают через $\mathcal{H}_h$.
\end{definition}

Очевидно, всякая гладкая выпуклая поверхность является 
гладким ежом, причём в этом случае 
введённые выше понятия опорной функции, 
главных радиусов кривизны и гауссовой кривизны
совпадают с общепринятыми.

Уже после того, как гладкие ежи завоевали право на существование
в геометрии, выяснилось, что статья \cite{LLR88} не была
первой, где они появились. 
В других обозначениях, в связи с другими задачами 
и только на плоскости они появлялись ещё в 1930-е годы.

Термин объясняется тем, что гладкие ежи параметризуются сферой
посредством своей нормали, т.\,е. каждой точке сферы
соответствует нормаль к ежу, причём  только одна.
Наверное, это есть идеал живого ежа (лат. Erinaceus europaeus):
его иглы должны быть направлены во все стороны 
(чтобы враг не прошёл), и в каждую сторону должно
быть направлено не более одной иглы (чтобы не наткнуться
на неё самому).

Пусть поверхность $S$ удовлетворяет условиям теоремы  \ref{1},
т.\,е. $S$ --- аналитическая выпуклая поверхность в $\mathbb{R}^3$,
гомеоморфная сфере,  причём существует постоянная $k$ такая, 
что для всех $\boldsymbol{x}\in S$ главные кривизны поверхности $S$
в точке $\boldsymbol{x}$ удовлетворяют неравенству
$k_1(\boldsymbol{x})\leqslant k\leqslant k_2(\boldsymbol{x})$ 
или, что то же самое,
$$(k_1(\boldsymbol{x})-k)(k_2(\boldsymbol{x})-k)\leqslant 0.
\eqno(3)$$
Обозначим опорную функцию поверхности $S$ через $f$,
а главные радиусы кривизны поверхности $S$ в точке $\boldsymbol{x}$
через $r_1(\boldsymbol{x})$ и $r_2(\boldsymbol{x})$.
Очевидно, $r_1(\boldsymbol{x})=1/k_1(\boldsymbol{x})$ 
и $r_2(\boldsymbol{x})=1/k_2(\boldsymbol{x})$.
 
Положим $r=1/k$, $h=f-r$. 
Прямое вычисление показывает, что главные радиусы кривизны 
гладкого ежа 
$x_h$ в точке $x_h(\boldsymbol{p})$
равны $r_1(x_f(\boldsymbol{p}))-r$ и $r_2(x_f(\boldsymbol{p}))-r$.
Следовательно, функция кривизны 
$R_h$ ежа $x_h$ задаётся формулой 
$(r_1(x_f(\boldsymbol{p}))-r)(r_2(x_f(\boldsymbol{p}))-r)$,
а значит в силу (3), неположительна всюду на сфере $\mathbb S^2$.
Тем самым теорема  \ref{1} вытекает из следующего утверждения:

\begin{theorem}\label{MM01}
Если функция кривизны $R_h$  
гладкого ежа $\mathcal{H}_h\subset\mathbb R^3$, 
соответствующего вещественно-аналитической функции 
$h:\mathbb S^2\to\mathbb R$, неположительна всюду на 
сфере $\mathbb S^2$, то $\mathcal{H}_h$ является точкой.
\end{theorem}

Легко понять, что и наоборот, теорема \ref{MM01} вытекает из
теоремы \ref{1}.
Независимое доказательство 
теоремы \ref{MM01} дано И. Мартинесом-Мором 
в статье \cite[теорема 3]{MM01}.

Доказательство И. Мартинеса-Мора 
основано на изучении 
полного прообраза $x_h^{-1}(\boldsymbol{s})$ 
специально подобранной точки 
$\boldsymbol{s}\in\mathcal{H}_h$.
Предполагая, что гладкий ёж $\mathcal{H}_h$ не является точкой, он
приходит к противоречию с тем фактом, что одномерное 
вещественное аналитическое множество имеет чётное число 
ветвей в каждой особой точке (см. \cite{Su71}).

Далее И. Мартинес-Мор отмечает, что 
его рассуждения не применимы к гладким ежам, 
имеющим особенности типа <<скрещенный колпак>>
(т.\,е. такие же особенности, как у конуса,
построенного над цифрой 8).
В силу упомянутого выше факта таких особенностей не
может быть у аналитических ежей.
Однако И. Мартинес-Мор явными формулами строит 
$C^2$-ежа с четырьмя особенностями типа 
<<скрещенный колпак>>
и всюду неположительной функцией кривизны 
$R_h=r_1r_2\leqslant 0$.
А именно, он показывает, что нужными свойствами
обладает поверхность $S\subset\mathbb{R}^3$,
являющаяся объединением графика функции
$(x^2-y^2)-\frac16(x^4-y^4)+\frac{1}{12}f(x,y)$
над множеством $|x|^{4/5}+|y|^{4/5}\leqslant 1$
и его зеркального отражения относительно плоскости
переменных $x,y$.
Здесь
$$
f(x,y)=\frac{xy\sqrt{(1-x^4-y^4)^{5/2}-25x^2y^2(x^8+y^8+3(x^4+y^4-
x^4y^4)+1)^{1/2}}}{1-x^4-y^4}.
$$

Прибавляя теперь к опорной функции гладкого ежа $S$
достаточно большое положительное число $R$
мы получим опорную функцию некоторой выпуклой поверхности $S_R$,
главные радиусы кривизны которой равны $r_1+R$ и $r_2+R$,
а значит, удовлетворяют неравенствам 
$r_1+R\leqslant R\leqslant r_2+R$.

Обратим внимание на то, что хотя функция
$f$ бесконечно дифференцируема в области 
$|x|^{4/5}+|y|^{4/5}< 1$,
ёж, получаемый из её графика указанным выше способом,
является лишь $C^2$-гладким.
Тем самым, вслед за И. Мартинесом-Мором, 
мы можем сказать, что выпуклая поверхность
$S_R$ является контрпримером к $C^2$-аналогу
теоремы \ref{1}. 

Однако, в рассуждениях И. Мартинеса-Мора,
приведённых в его статье \cite{MM01},
есть серьёзный изъян: неотрицательность функции кривизны 
$R_h$ ежа $S$ проверена только численно.

Таким образом в 2001 году сложилась следующая
парадоксальная ситуация: с одной стороны, 
в работе \cite{MM01} 2001 года
И. Мартинес-Мор показал (хотя и не <<чисто>>),
что  $C^2$-аналог теоремы \ref{1} неверен, а 
с другой стороны, А.В. Погорелов (1919--2002) 
в 1998 году анонсировал 
(хотя и не привёл подробного доказательства)
в <<Докладах РАН>> \cite{Po98}, 
что аналог теоремы \ref{1} справедлив для
$C^2$-гладких поверхностей.
Тем самым, в 2001 году стало ясно, что как минимум
одна из работ \cite{Po98} и \cite{MM01} ошибочна.
Весь вопрос в том --- какая. 
В 2005 году на этот вопрос ответила Г.Ю. Панина
в статье \cite{Pa05} (а затем и в докторской 
диссертации \cite{Pa06b}).
Мы расскажем о её работе в следующем разделе.

\section{Многогранные ежи}

До сих пор мы имели дело только с гладкими ежами,
т.\,е. с огибающими семейств плоскостей
$
\{\boldsymbol{x}\in\mathbb R^n\vert 
(\boldsymbol{x},\boldsymbol{p})=h(\boldsymbol{p})\}
$
для заданного гладкого отображения 
$h:\mathbb S^{n-1}\to\mathbb R$.
Теперь обсудим многогранный аналог этого понятия,
который будем называть многогранным ежом.
Нужно только иметь ввиду, что
разные авторы определяют этот объект немного по-разному,
см., например, \cite{LLR88, Al04, Pa06a}. 

Мы начнём с наиболее наглядного определения 
многогранного ежа,
заимствованного из \cite{Al04}.

\begin{definition}
Многогранным ежом в $\mathbb R^3$ называется многогранная поверхность $P$ в
$\mathbb R^3$, каждая грань которой оснащена единичным вектором 
$\boldsymbol{n}_j$, перпендикулярным
к $Q_j$, таким образом, что если грани $Q_j$ и $Q_k$ имеют общее ребро, то
$\boldsymbol{n}_j + \boldsymbol{n}_k \neq 0$, а значит концы векторов 
$\boldsymbol{n}_j$ и $\boldsymbol{n}_k$ могут быть соединены
на единичной сфере $\mathbb S^2\subset\mathbb R^3$ единственной кратчайшей.
Обозначив эту кратчайшую через $l_{jk}$, потребуем дополнительно чтобы при любых
значениях индексов $j$ и $k$, соответствующих смежным граням, две дуги вида $l_{jk}$
либо не пересекались, либо имели один общий конец.
Потребуем, наконец, чтобы объединение всех дуг $l_{jk}$ образовывало разбиение сферы
$\mathbb S^2$ на многоугольники, хотябы и не выпуклые.
\end{definition}

Несмотря на громоздкость этого определения, 
многогранные ежи служат довольно естественным
обобщением выпуклых многогранников.
В частности, граница всякого выпуклого многогранника 
в $\mathbb{R}^3$, оснащённая, скажем,
внешними нормалями к 2-мерным граням, 
является многогранным ежом.

\begin{figure}\center           
 \includegraphics{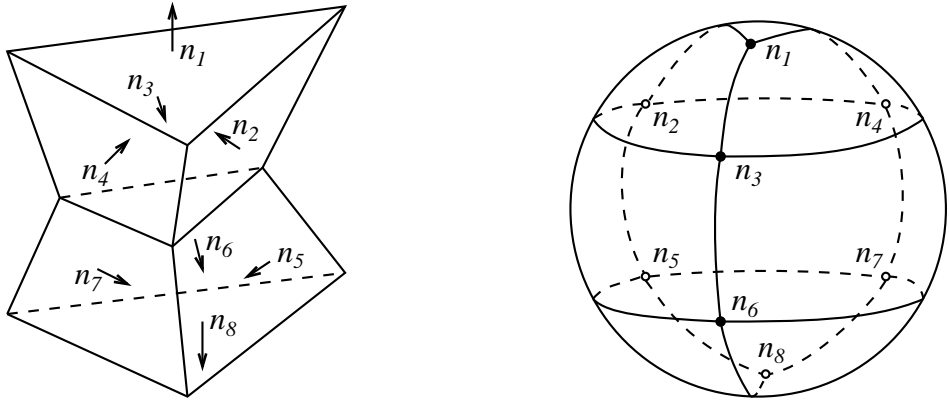}
 \caption{Невыпуклый многогранный ёж в $\mathbb R^3$ 
 и соответствующее разбиение сферы на многоугольники}
\end{figure}

Приведём пример невыпуклого многогранного ежа.
Для этого рассмотрим правильный тетраэдр $\Delta$ в $\mathbb{R}^3$, отсечём от $\Delta$ меньший 
тетраэдр плоскостью $\tau$, параллельной одной из граней $\Delta$, и рассмотрим 
объединение оставшегося усечённого тетраэдра с его симметричным образом 
относительно плоскости $\tau$. 
Невыпуклую многогранную поверхность, являющуюся границей этого объединения, 
обозначим через $P$ (см. рис. 2). 
Оснастив грани $P$, параллельные $\tau$, внешними нормалями, а все остальные грани
--- внутренними, мы превратим $P$ в многогранного ежа.
На рис. 2 изображено также разбиение сферы $\mathbb S^2$, соответствующее этому ежу.

В \cite{Al04}, в частности, доказано, что на 
определённые выше многогранные 
ежи переносится классическая теорема Минковского
о существовании и единственности выпуклого многогранника
с заданными направлениями и площадями граней.

Оказывается, что приведённое выше определение многогранного ежа
является лишь одним из многих равносильных 
(или почти равносильных) описаний этого объекта. 
Подобные объекты более ли менее независимо возникали в 
разных областях математики. 
После работы А.В. Пухликова и А.Г. Хованского \cite{PKh92} 
их обычно называют виртуальными многогранниками.  
Следуя \cite{Pa06b} укажем для них 
следующие равносильные представления:

(1) Виртуальный многогранник это элемент группы Гротендика
полугруппы выпуклых многогранников (где групповой операцией 
служит солжение по Минковскому $\otimes$) в евклидовом
пространстве, т.\,е. формальное выражение вида
$K\otimes M^{-1}$, где $K$ и $M$ --- выпуклые многогранники,
детали см. в  \cite{PKh92}.

(2) Виртуальный многогранник есть многогранная функция,
т.\,е. конечная линейная комбинация с целочисленными коэффициентами
характеристических функций выпуклых многогранников,
детали см. в  \cite{PKh92}.

(3) Виртуальный многогранник есть кусочно-линейная 
положительно-однородная функция, заданная в евклидовом пространстве,
детали см. в  \cite{PKh92}.

(4) Виртуальный многогранник это пара $(F,\Sigma)$,
где $F$ --- некоторая замкнутая многогранная поверхность
(возможно с самопересечениями и самоналожениями)
с коориентированными гранями, а $\Sigma$ ---
ассоциированный с $F$ веер,
детали см. в  \cite{Pa02, Pa05, Pa06a}.

В \cite{Pa06b} показано, что с каждым виртуальным многогранником
можно естественным образом ассоциировать двойственный объект ---
график его опорной функции. 
Существенным вкладом Г.Ю. Паниной в обсуждаемую тематику
стало то, что в своих работах \cite{Pa02, Pa05, Pa06a, Pa06b}
она ввела понятие гиперболического виртуального
многогранника (т.\,е. такого виртуального многогранника, для
которого график сужения его опорной функции на любую плоскость
является седловой поверхностью, т.\,е. поверхностью, от которой 
никакая плоскость не отсекает компактного куска).
Далее, она научилась строить разнообразные примеры 
гиперболических
виртуальных многогранников.
И наконец, доказала общую теорему о сглаживании 
гиперболических виртуальных многогранников.
Тем самым она научилась, отправляясь от данного
гиперболического виртуального многогранника,
строить гладкого ежа.

В результате Г.Ю. Панина смогла построить целую серию
$C^\infty$-гладких ежей, каждый из которых порождает контрпример
к $C^\infty$-аналогу теоремы \ref{1} ---
нужно лишь сложить (по Минковскому) этот ёж и сферу
достаточно большого радиуса, как это и делал И. Мартинес-Мор.
Тем самым, Г.Ю. Панина смогла не только подтвердить 
правильность теоремы И. Мартинеса-Мора, но и 
сумела распространить её на $C^\infty$-поверхности.

Более того, Г.Ю. Панина явно указала на ошибку в рассуждениях
А.В. Погорелова \cite{Po98}, пропустившего один из случаев в своих
рассуждениях. Как оказалось, именно этот случай реализуется 
в примерах, построенных и И. Мартинес-Мором и Г.Ю. Паниной.

\section{Заключение}

Выше мы постарались показать, что 
классическая теорема А.Д. Александрова \ref{1}
связана со многими вопросами теории поверхностей,
выпуклой геометрии, дифференциальных уравнений, алгебры и 
даже алгебраической геометрии. 
Очевидно, тут есть и много открытых вопросов.
В самом общем виде программа исследований может быть 
сформулирована так:
до каких пределов теория выпуклых тел переносится на ежи?

Остаётся надеяться, что со временем
мы будем узнавать всё новые и новые результаты, полученные
в этом направлении.

\bigskip

\section*{Приложение: {\small перевод с французского статьи}}
\centerline{{\bf Контрпример к одной гипотетической характеризации
сферы}\footnote{Оригинал этой статьи вышел на французском языке:
Yves Martinez-Maure, {\it Contre-exemple \`{a} une caract\'{e}risation
conjectur\'{e}e de la sphe\`{e}re} // C. R. Acad. Sci. Paris, S\'{e}rie~I, 
{\bf 332}:1 (2001) 41--44. Copyright \copyright \,2001 Acad\'{e}mie 
des sciences.
Издательство Elsevier Masson SAS. Все права сохранены.
В.А. Александров благодарит издательство Эльзевир Массон
за предоставленную возможность опубликовать этот перевод.}}

\centerline{{\bf Ив Мартинес-Мор}}

\small{
\centerline{1, rue Auguste-Perret, 92500 Rueil-Malmaison, France}

\centerline{(Статья получена 19 октября 2000, принята к печати 23 октября 2000)}

\centerline{Статья представлена Марселем Берже}

\bigskip

\noindent{\bf Аннотация}

\noindent{Уже} давно была высказана гипотеза, что замкнутая выпуклая поверхность класса 
$C^2_+$, главные кривизны $K_1$, $K_2$ которой удовлетворяют неравенству
$(K_1-c)(K_2-c)\leqslant 0$ с некоторой постоянной $c$, обязана быть сферой.
Частные результаты были получены А.Д. Александровым, 
Х.Ф. Мюнцнером и Д. Коутроуфиотисом.\newline
Мы переформулируем эту гипотезу в терминах ежей и 
строим контрпример.
Кроме того, мы доказываем эту гипотезу для поверхностей 
постоянной ширины и даём 
новое доказательство для аналитических поверхностей.

\medskip

\noindent{\bf 1. Введение}

Уже давно была высказана гипотеза, что любая  выпуклая замкнутая поверхность класса 
$C^2_+$ (т.\,е. класса $C^2$ с гауссовой кривизной $>0$), главные кривизны $K_1$ и 
$K_2$ которой удовлетворяют неравенству $(K_1-c)(K_2-c)\leqslant 0$ с некоторой 
постоянной $c$, с необходимостью является сферой.
Для аналитических поверхностей это было установлено А.Д. Александровым [1п, 2п] 
и Х.Ф. Мюнцнером [7п].
Эта гипотеза также проверена для поверхностей вращения [6п] и для более общего 
класса поверхностей, у которых одна из ортогональных проекций является окружностью [3п].

В настоящей заметке мы переформулируем эту гипотезу в терминах ежей и строим контрпример.
Понятие ежа служит обобщением понятия выпуклого тела класса $C^2_+$.
Каждому выпуклому телу $K\subset\mathbb R^{n+1}$ сопоставляется опорная функция
$h_K:\mathbb S^n\to\mathbb R$, 
$p\mapsto h_K(p)=\mbox{max\, }\{ \langle m,p\rangle \ | \ m\in K\}$.
Когда $K$ принадлежит классу $C^2_+$, $h_K$ принадлежит $C^2$ и определяет границу $K$
как огибающую семейства гиперплоскостей, заданных уравнениями
$\langle x,p\rangle =h_K(p)$.
Даже если функция $h\in C^2(\mathbb S^n;\mathbb R)$ не обязательно является опорной функцией
некоторого выпуклого тела, с ней всё равно можно связать огибающую $\mathcal H_h$
семейства гиперплоскостей, заданных уравнением $\langle x,p\rangle =h(p)$.
Эта огибающая параметризована с помощью $x_h:\mathbb S^n\to\mathbb R^{n+1}$,
$p\mapsto (\mbox{grad\, } h)(p)+h(p)p$.
Отображение $x_h$ может быть интерпретировано как обратное к отображению Гаусса для
$\mathcal H_h$:  каждый регулярный кусок $\mathcal H_h$ допускает введение ориентации
с помощью нормали таким образом, что $p$ оказывается вектором нормали в точке 
$x_h(p)$.
В работе Р. Ланжевена, Г. Левита и Г. Розенберга [4п] гиперповерхность $\mathcal H_h$
названа ежом опорной функции $h$.

Пусть $S$ есть некоторая замкнутая выпуклая поверхность класса $C^2_+$,
главные кривизны  $K_1$ и  $K_2$ которой удовлетворяют неравенству 
$(K_1-c)(K_2-c)\leqslant 0$ при некоторой постоянной $c$.
Положив $r=c^{-1}$, мы можем переписать это неравенство в виде
$(R_1-r)(R_2-r)\leqslant 0$, где $R_1$ и $R_2$ --- главные 
радиусы кривизны поверхности $S$.
Обозначим через $f$ опорную функцию поверхности $S$.
Так как $x_f:\mathbb S^2\to S$ является обратным к гауссову отображению поверхности $S$,
то $R_1(p)$ и $R_2(p)$ являются собственными значениями его касательного отображения
$T_px_f$ в точке $p$.
Положив $h=f-r$, выводим, что $R_1(p)-r$ и $R_2(p)-r$ являются собственными значениями 
касательного отображения  к $x_h$ в точке $p$.
Исходная гипотеза таким образом может быть переформулирована так: <<Если функция
кривизны $R_h(p)=\mbox{det \,} (T_px_h)$ ежа $\mathcal H_h\subset\mathbb R^3$
неположительна на сфере $\mathbb S^2$, то $\mathcal H_h$ с необходимостью сводится 
к точке>>.

Изучение проекций ежа $\mathcal H_h$ даёт нам информацию о сферическом представлении 
сингулярностей.
Мы получим несколько известных частных результатов и новое доказательство для 
поверхностей постоянной ширины.
В общем случае метод не проходит из-за существования сингулярностей типа <<скрещенный колпак>>.
Мы предъявляем ежа с особенностью типа <<скрещенный колпак>> и функцией кривизны $\leqslant 0$ 
(моделируемый на $\mathbb S^2$ без некоторых полуокружностей).
Наконец, мы строим нетривиальный ёж с функцией кривизны $\leqslant 0$ на $\mathbb S^2$.
В силу проективной двойственности это влечёт существование не вполне геодезической 
2-сферы, погруженной в $\mathbb S^3$ с внешней кривизной $\leqslant 0$.

\medskip

\noindent{\bf 2. Изучение специальных случаев}

Пусть $\mathcal H_h\subset\mathbb R^{n+1}$ --- ёж и 
$x\in\mathbb R^{n+1}\setminus\mathcal H_h$.
Индексом $i_h(x)$ называется алгебраическое число пересечений ориентированной 
полупрямой, выходящей из $x$, с ежом $\mathcal H_h$, снабжённым ориентацией
(это число не зависит от конкретной полупрямой из некоторого открытого плотного
множества направлений).
Функция кривизны ежа $\mathcal H_h\subset\mathbb R^3$ связана с индеком плоского ежа,
получаемого из $\mathcal H_h$ посредством ортогонального проектирования.
Для всех $n\in\mathbb S^2$ ограничение функции $h$ на окружность
$\mathbb S^1_n=\mathbb S^2\cap n^\perp$ является опорной функцией $h_n$
<<проекции>> ежа $\mathcal H_h$ на $n^\perp$: 
$\mathcal H_{h_n}=(\pi_n\circ x_h)(\mathbb S^1_n)$,
где $\pi_n$ является ортогональной проекцией на $n^\perp$.
Для индекса этой проеции справедлива следующая 

ТЕОРЕМА 1 ([5п]). {\it Для каждого регулярного значения $x$ ограничения
$\pi_n\circ x_h$ на полусферу 
$\mathbb S^2_n=\{p\in\mathbb S^2\ | \ \langle p,n\rangle\geqslant 0\}$
верно равенство $i_{h_n}(x)=\nu^n_h(x)^+-\nu^n_h(x)^-$,
где $\nu^n_h(x)^+$ (соотв. $\nu^n_h(x)^-$) есть число точек 
$p\in\mathbb S^2_n$ таких, что $x_h(p)$ является эллиптической (соотв. гиперболической)
точкой $\mathcal H_h$, лежащей на прямой $\{ x\} +\mathbb R n$.}

Следовательно, если $\mathcal H_h$ имеет функцию кривизны $\leqslant 0$ и 
не сводится к одной точке, то $i_{h_n}\leqslant 0$ и не является тождественным нулём.
Так как $i_{h_n}\geqslant 0$ если $\mathcal H_{h_n}$ является окружностью или точкой,
то это доказывает гипотезу для поверхностей, хоть одна из проекций которых
является окружностью.
Введём понятие точки ежа, которым нам предстоит пользоваться ниже.
Каждому ежу $\mathcal H_h\subset\mathbb R^{n+1}$ соответствует псевдометрика,
определённая на $\mathbb S^n$ равенством 
$\rho_h(p,q)=\mbox{inf \,} (\{ L(x_h\circ \gamma) \ | \ \gamma\in C_{p,q}\}$,
где $C_{p,q}$ есть множество кусочно $C^1$-гладких кривых  с концами $p$ и $q$,
а $L(x_h\circ\gamma)$ является длиной кривой $x_h\circ\gamma$.
Точкой ежа $\mathcal H_h$ по определению является точка пространства 
$(\mathbb S^n/\sim )=\{ [p] \ | \  p\in\mathbb S^n\}$,
где $p\sim q$ если и только если $\rho_h (p,q)=0$,
снабженном метрикой $d_h([p],[q])=\rho_h(p,q)$.
Но для простоты мы не различаем точку $[p]$ и её геометрическую реализацию $x_h(p)$.

СЛЕДСТВИЕ 1. {\it Если $\mathcal H_h\subset\mathbb R^3$ имеет функцию кривизны $\leqslant 0$,
то (i) для всех $x\in n^\perp -\mathcal H_{h_n}$ существует ровно $-2i_{h_n}(x)$
точек $\mathcal H_h$, геометрические реализации которых лежат на прямой $\{ x\} +\mathbb R n$;
(ii) для каждой точки $x\in\mathcal H_{h_n}$ каждая точка ежа $\mathcal H_h$,
геометрическая реализация которой лежит на прямой $\{ x\} +\mathbb R n$,
является классом некоторой точки окружности $\mathbb S^1_n=\mathbb S^2\cap n^\perp$.}

Пусть $\mathcal H_h$ является ежом в $\mathbb R^{n+1}$.
Говорят, что $s\in x_h(\mathbb S^n)$ является экстремальной точкой для $\mathcal H_h$
 в направлении $n\in\mathbb S^n$, если 
$\langle x_h(p),n\rangle\leqslant \langle s,n\rangle$ для всех $p\in\mathbb S^n$.

ЛЕММА 1. {\it Пусть $s$ является экстремальной точкой ежа $\mathcal H_h\subset\mathbb R^2$
в направлении $n\in\mathbb S^1$.
Тогда на $\mathbb S^1$ компоненты связности 
множества $x^{-1}_h(s)$ отделяются одна от другой 
точками множества $\{ n,-n\}$.}

{\it Доказательство.}
Без ограничения общности можно считать, что $s=(0,0)$ и $n=(0,1)$.
Тогда опорная прямая с нормалью $u_\theta =(\cos\theta, \sin\theta )\notin \{ -n,n\}$
пересекает $n^\perp$ в точке $(x(\theta), 0)$, где
$x(\theta)=h(n_\theta )/\cos\theta$.
Если внутренность дуги $\Gamma\subset\mathbb S^1$ с концами $\alpha ,\beta\in x^{-1}_h(s)$
не пересекает $\{ -n,n\}$, то $\Gamma\subset x^{-1}_h(s)$.
В самом деле, тогда мы бы имели $x(\theta)\to 0$ если $u_\theta\to\alpha$ или $\beta$,
$x'(\theta)=\langle x_h(u_\theta),n\rangle/\cos^2\theta\leqslant 0$ если 
$u_\theta\in\Gamma-\{\alpha,\beta\}$ и, следовательно, $h=0$ на $\Gamma$.
$\square$

ТЕОРЕМА 2. {\it Пусть ёж $\mathcal H_h\subset\mathbb R^3$ имеет функцию кривизны $\leqslant 0$.
Если $\mathcal H_h$ не является точкой, то $\mathcal H_h$ допускает в некотором напрвлении
$n\in\mathbb S^2$ экстремальную точку $s$, сферическое представление $x^{-1}_h(s)$
которой не пересекает одну из полуокружностей, соединяющих $-n$ и $n$ на $\mathbb S^2$.}

{\it Доказательство.}
Если $\mathcal H_h$ не является точкой, то существует направление $n\in\mathbb S^2$
в котором $\mathcal H_h$ имеет две различные экстремальные точки $s_1$ и $s_2$.
Положим $p=(s_2-s_1)/\| s_2-s_1\| $ и рассмотрим проекцию $\mathcal H_h$ на $p^\perp$.
Общая проекция $s$ точек $s_1$ и $s_2$ является экстремальной точкой $\mathcal H_{h_p}$
в направлении $n$.
Одна из компонент связности множества $x^{-1}_{h_p}(s)$ не пересекает сразу
(одновременно) $x^{-1}_h(s_1)$ и $x^{-1}_h(s_2)$ так как $\mathcal H_h$
не содержит сегмента $[s_1,s_2]$.
Однако $x^{-1}_h(s_1)$ и $x^{-1}_h(s_2)$ пересекают $x^{-1}_{h_p}(s)$ 
(следствие~1) и их компоненты разделяются множеством $\{ -n,n\}$ на $\mathbb S^1_p$
(лемма~1), следовательно, $x^{-1}_h(s_1)$ (соотв. $x^{-1}_h(s_2)$) не пересекает 
одну из полуокружностей, соединяющих $-n$ и $n$ на $\mathbb S^1_p$. $\square$

Ёж $\mathcal H_h\subset\mathbb R^3$ называется аналитическим (соотв. проективным)
если его опорная функция $h$ является аналитической (соотв. антисимметричной:
для каждого $p\in\mathbb S^n$, $h(-p)=-h(p)$).

ТЕОРЕМА 3. {\it Если $\mathcal H_h\subset\mathbb R^3$ является аналитическим
(соотв. проективным) ежом с функцией кривизны $\leqslant 0$, то $\mathcal H_h$
является точкой.}

{\it Доказательство.}
Пусть $s$ является экстремальной точкой ежа $\mathcal H_h$ в напрвлении 
$n\in\mathbb S^2$ и $[m]$ является точкой ежа $\mathcal H_h$, для которой
$s$ является геометрической реализацией.
Согласно теореме 2 достаточно проверить, что $x^{-1}_h(s)$ пересекает
каждую полуокружность, соединяющую $-n$ и $n$ на $\mathbb S^2$.
Согласно следствию 1, $[m]$ пересекает каждую большую окружность, проходящую
через $n$.
Если $\mathcal H_h$ является проективным ежом, то $x^{-1}_h(s)$ инвариантно
относительно антиподального преобразования и результат достигнут.
Допустим $\mathcal H_h$ является аналитическим ежом, не сводящимся к точке.
Поскольку отображение $\phi : \mathbb S^2\to\mathbb R$, $u\mapsto\langle x_h(u), n\rangle$
аналитично, компоненты связности $\phi (u)=\langle s,n\rangle$ являются либо
точками, либо простыми замкнутыми дугами, либо погружениями некоторых связных графов,
все вершины которых имеют чётную степень [8п].
Обозначим через $C$ ту компоненту связности, которая содержит $[m]$.
Поскольку $[m]$ пересекает все большие дуги, проходящие через $n$, и
пересекает самое большее один раз каждую большую полуокружность, соединяющую
$-n$ и $n$ (лемма 1), то либо $C$ пересекает $\{ -n,n\}$, либо $C$ является
простой дугой, разделяющей $-n$ и $n$ на $\mathbb S^2$.
Однако, $C\subset x^{-1}_h(s)$ (иначе некторая проекция ежа $\mathcal H_h$
содержала бы нетривиальный сегмент).
Следовательно, $x^{-1}_h(s)$ пересекает каждую полуокружность, соединяющую
$-n$ и $n$ на $\mathbb S^2$, что противоречит теореме 2. $\square$

СЛЕДСТВИЕ 2. {\it Гипотеза верна для аналитических поверхностей (соотв. для
поверхностей постоянной ширины).}

{\it Доказательство.}
Случай аналитических поверхностей непосредственно вытекает из теоремы 3.
Если $S$ имеет постоянную ширину, то $h=f-r$ имеет вид $g+k$, где $g$
антисимметрично, а $k$ постоянно.
Тогда мы имеем $R_h=R_g+2kR_{(1,g)}+k^2$, где $R_{(1,g)}$ обозначает
среднее арифметическое главных радиусов кривизны ежа $\mathcal H_g$.
Отсюда замечая, что неравенство $R_h\leqslant 0$, примененное к антиподальным особым
точкам $x_g$, даёт $k^2\leqslant 0$, получаем, что $\mathcal H_h$ является проективным ежом.
$\square$

\medskip

\noindent{\bf 3. Контрпример}

Ни метод параграфа 2, ни предшествующие работы не позволяют принимать во внимание
особенности типа скрещенного колпака.
Отображение $X:\mathbb R^2\to\mathbb R^3$, $(u,v)\mapsto (x,y,z)=r^4(u,1,uv)$, 
где $r=\sqrt{u^2+v^2}$, определяет скрещенный колпак, который можно рассматривать как
ёж $\mathcal H_h$ (моделируемый на сфере $\mathbb S^2$ с удаленной полуокружностью
$z=0$, $x\leqslant 0$) с функцией кривизны $\leqslant 0$.
На $\mathbb S^2$ множество критических точек ежа $\mathcal H_h$ является полуокружностью 
$y=0$, $x\geqslant 0$.
Являясь сферическим представлением экстремальной точки в направлении $n=(0,-1,0)$,
эта полуокружность не пересекает геодезические линии сферы $\mathbb S^2$, которые
соединяют $-n$ и $n$ и расположены в полусфере $x<0$.
Метод параграфа 2 не может быть, следовательно, применён к изучению таких 
сингулярностей.
Регулярная часть ежа $\mathcal H_h\subset\mathbb R^3$ является частью алгебраической 
поверхности, заданной уравнением
$x^5y^5-(x^4+y^2z^2)^2=0$,
содержащейся в полупространстве $y>0$.
Следовательно, она может быть получена склеиванием графика функции
$f(x,y)=(x/y)\sqrt{y^{5/2}-x^2}$ с его симметричным образом относительно плоскости
$z=0$.
Проверка не представляет особых трудностей.
Этот пример показывает, что неаналитический ёж может иметь аналитическую 
параметризацию.

Приведенный пример наводит на мысль соединить 4 скрещенных колпака так, чтобы
получить ёж, с функцией кривизны $\leqslant 0$ на $\mathbb S^2$.
Начнём со склеивания графика функции
$$
f(x,y)=\frac{xy}{1-x^4-y^4}\sqrt{(1-x^4-y^4)^{5/2}-25x^2y^2(x^8+y^8+3(x^4+y^4-x^4y^4)+1)^{1/2}},
$$
где $(x,y)\in D=\{ (u,v)\in\mathbb R^2 | |u|^{4/5} +|v|^{4/5} \leqslant 1\}$,
с его симметричным образом относительно плоскости $z=0$
(заметим, что подкоренное выражение в точности зануляется на границе области $D$).
Полученная поверхность образована четырьмя скрещенными колпаками, но имеет кривизну $\geqslant 0$.
Наконец, чтобы исключить кривизну $\geqslant 0$, добавим к $f$ функцию вида 
$g(x,y)=a(x^2-y^2)+b(x^4-y^4)$ с $a\neq 0$ и $a+6b=0$, 
для которых график функции $g:D\to\mathbb R$ имеет кривизну $<0$ всюду, кроме особых точек на 
границе.
Иными словами, рассмотрим параметрическое семейство поверхностей
$(\mathcal S_t)_{t\in\mathbb R^*_+}$ , определяемое формулами
$X_t:D\times\{-1,1\}\to\mathbb R^3$, $(x,y,\varepsilon)\mapsto (x,y,(x^2-y^2)-\frac{1}{6}(x^4-y^4)+t\varepsilon f(x,y))$.

Вычисления на компьютере показывают, что существует интервал значений параметра $t$, 
для которых $\mathcal S_t$ является поверхностью с гауссовой кривизной $<0$.
Среди этих значений мы выбираем $t=1/12$.
Поверхность $\mathcal S_{1/12}$ является частью алгебраической поверхности.
Она симметрична относительно плоскостей $x=0$ и $y=0$ и относительно прямых
$x=y=0$ и $z=0$, $x=y$ (соотв. $x=-y$).
Поскольку она является склейкой двух графиков, расположенных над $D$ 
(чьи границы являются ежами с опорной функцией $(u,v)\mapsto uv(u^4+v^4)^{-1/4}$,
$(u,v)\in\mathbb S^1\subset\mathbb R^2$) строгая отрицательность её гауссовой кривизны
(и, следовательно, локальная инъективность её гауссова отображения) влечёт то, что
она может рассматриваться как некоторый ёж $\mathcal H_h$.
Изучим его класс дифференцируемости.
Вычисления на компьютере показывают, что множество особых точек ежа $\mathcal H_h$
образовано на $\mathbb S^2$ четырьмя полуокружностями:
(i) $3x+4z=0$, $y\geqslant 0$;
(ii) $3y-4z=0$, $x\geqslant 0$;
(iii) $3x-4z=0$, $y\leqslant 0$;
(iv) $3y+4z=0$, $x\leqslant 0$.
Сразу проверяется, что $h$ принадлежит $C^\infty$ вне сингулярного множества
и принадлежит $C^1$ на всём $\mathbb S^2$ (заметим, что $x_h$ является ограничением на $S^2$
градиента функции $\varphi (u)=\| u\| h(u/\| u\|)$, $u\in\mathbb R^3-\{0_{\mathbb R^3}\})$.
Можно доказать также, что $h$ принадлежит $C^2$ на $\mathbb S^2$ и заметить, что на $\mathbb S^2$
собственные значения отображения $(\mbox{hess\,} \varphi )(p)$ (0 и главные радиусы кривизны
ежа $\mathcal H_h$  в точке $x_h(p)$) стремятся к нулю при подходе к сингулярному множеству
(вычисления показывают, что $R_h$ и $R_{(1,h)}$ стремятся к нулю при приближении
к сингулярному множеству).

\medskip

\noindent{\bf Цитированная литература}

\noindent{[1п]}  А.Д. Александров,
{\it О теоремах единственности для замкнутых поверхностей},
Докл. АН СССР {\bf 22}, \No 3 (1939), 99--102.

\noindent{[2п]} А.Д. Александров, 
{\it О кривизне поверхностей}, 
Вестн. ЛГУ. Сер. математики, механики, астрономии. 
\No 19, вып. 4 (1966), 5--11.

\noindent{[3п]} D. Koutroufiotis,
{\it On a conjectured characterization of the sphere},
Math. Ann. {\bf 205} (1973), 211--217.

\noindent{[4п]} 
R. Langevin, G. Levitt, H. Rosenberg,
{\it H\'{e}rissones et multih\'{e}rissons 
(enveloppes parametre\'{e}s par leur application de Gauss)}, 
Singularities, Banach Center publ. {\bf 20} (1988), 245--253.

\noindent{[5п]} 
Y. Martinez-Maure,
{\it Indice d'un h\'{e}risson: \'{E}tude et applications},
Publ. Mat., Barc. {\bf 44}, no. 1 (2000), 237--255.

\noindent{[6п]} 
H.F. M\"{u}nzner,
{\it \"{U}ber eine spezielle Klasse von Nabelpunkten und analoge Sin\-gu\-la\-ri\-t\"{a}\-ten 
in der zentroaffinen Fl\"{a}chentheorie},
Comment. Math. Helv. {\bf 41} (1966), 88--104. 

\noindent{[7п]}
H.F.  M\"{u}nzner,
{\it \"{U}ber Fl\"{a}chen mit einer Weingartenschen Ungleichung},
Math. Z. {\bf 97} (1967), 123--139.

\noindent{[8п]}
D. Sullivan,
{\it Combinatorial invariants of analytic spaces},
В кн.: Proc. Liverpool Singulari\-ties-Sympos. I. 
Lect. Notes Math. {\bf 192} (1971), 165--168.
}%end of small


\bigskip
\begin{thebibliography}{99}

\bibitem{AR04} U. Abresch, H. Rosenberg, 
{\it A Hopf differential for constant mean curvature surfaces in 
$\mathbf{S}^2\times\mathbf{R}$ and 
$\mathbf{H}^2\times\mathbf{R}$,}
Acta Math. {\bf 193}, no. 2 (2004), 141--174.
doi: 10.1007/BF02392562

\bibitem{AD38} А.Д. Александров,
{\it Одна общая теорема единственности для замкнутых поверхностей},
Докл. АН СССР. {\bf 19}, № 4 (1938), 233--236. 
(Переиздание в книге: А.Д. Александров, 
{\it Геометрия и приложения}, Новосибирск: Наука, 2006
(Избранные труды, Т. 1. ISBN 5-02-032427-2).) 

\bibitem{AD66} А.Д. Александров,
{\it О кривизне поверхностей},
Вестн. ЛГУ. № 19. Сер. математики, механики и астрономии. Вып. 4 (1966), 5--11.

\bibitem{Al04} V. Alexandrov, 
{\it Minkowski-type and Alexandrov-type theorems for polyhedral herissons},
Geom. Dedicata {\bf 107} (2004), 169--186. 
doi: 10.1007/s10711-004-4090-3

\bibitem{BVK73} И.Я. Бакельман, А.Л. Вернер, Б.Е. Кантор,
{\it Введение в дифференциальную геометрию <<в целом>>},
М.: Наука, 1973.

\bibitem{BF34} T. Bonnesen, W. Fenchel,
{\it Theorie der konvexen K\"{o}rper,}
Springer, Berlin (1934).
(Русский перевод: Т. Боннезен, В. Фенхель,
{\it Теория выпуклых тел}, М.: Фазис, 2002.
ISBN: 978-5-7036-0075-7)

\bibitem{GM10} G. Ganchev, V. Mihova,
{\it On the invariant theory of Weingarten 
surfaces in Euclidean space,}
J. Phys. A, Math. Theor. {\bf 43}, no. 40 (2010), 
Article ID 405210, 27 p. 
doi: 10.1088/1751-8113/43/40/405210 

\bibitem{GM03} P. Guan, X.-N. Ma,             
{\it The Christoffel--Minkowski problem. I: 
Convexity of solutions of a Hessian equation},
Invent. Math. {\bf 151}, no. 3 (2003), 553--577.
doi: 10.1007/s00222-002-0259-2

\bibitem{LLR88} 
R. Langevin, G. Levitt, H. Rosenberg,
{\it H\'{e}rissones et multih\'{e}rissons 
$($enveloppes parametre\'{e}s par leur application de Gauss}),
Singularities, Banach Center Publ. 
{\bf 20} (1988), 245--253.

\bibitem{Li00} H. Liebmann, 
{\it Ueber die Verbiegung der geschlossenen 
Fl\"{a}chen positiver Kr\"{u}mmung},
Math. Ann. {\bf 53}, no. 1 (1900), 81--112.

\bibitem{MM01} Y. Martinez-Maure,             
{\it Contre-exemple \`{a} une caract\'{e}risation conjectur\'{e}e de la sph\`{e}re}, 
C. R. Acad. Sci., Paris, S\'{e}r. I, Math. {\bf 332}, no. 1 (2001), 41--44.
doi: 10.1016/S0764-4442(00)01756-0
(Русский перевод дан в <<Приложении>> к настоящей статье.)

\bibitem{MM10} Y. Martinez-Maure,             
{\it New notion of index for hedgehogs of $\mathbb{R}^3$  
and applications}, 
Eur. J. Comb. {\bf 31}, no. 4 (2010), 1037--1049.
doi: 10.1016/j.ejc.2009.11.015

\bibitem{Mu67} H.F. M\"{u}nzner,
{\it \"{U}ber Fl\"{a}chen mit einer Weingartenschen Ungleichung,}
Math. Zeitschrift {\bf 97}, no. 2 (1967), 123--139.    
{\it Bemerkung zur Arbeit ,,\"{U}ber Fl\"{a}chen mit einer Weingartenschen Ungleichung``}, ibid {\bf 100}, no. 5 (1967), 416.

\bibitem{Pa02} Г.Ю. Панина,
{\it Виртуальные многогранники и классические вопросы геометрии},
Алгебра и анализ {\bf 14}, № 5 (2002), 152--170.

\bibitem{Pa05} G. Panina, 
{\it New counterexamples to A.D. Alexandrov's hypothesis},
Adv. Geom. {\bf 5}, no. 2 (2005), 301--317.
doi: 10.1515/advg.2005.5.2.301

\bibitem{Pa06a} G. Panina, 
{\it On hyperbolic virtual polytopes and hyperbolic fans},
Cent. Eur. J. Math. {\bf 4}, no. 2 (2006), 270--293.
doi: 10.2478/s11533-006-0006-9

\bibitem{Pa06b} Г.Ю. Панина,
{\it Виртуальные многогранники},
Дисс.\dots д-ра физ.-мат. наук, Ст.-Петербург, 2006. 

\bibitem{Po98} А.В. Погорелов,
{\it Решение проблемы А.Д. Александрова},
Докл. РАН {\bf 360}, № 3 (1998), 317--319.

\bibitem{PKh92} А.В. Пухликов, А.Г. Хованский,
{\it Конечно-аддитивные меры виртуальных многогранников},
Алгебра и анализ {\bf 4}, № 2 (1992), 161--185.

\bibitem{Su71} D. Sullivan,
{\it  Combinatorial invariants of analytic spaces},
In the book: {\it Proceedings of Liverpool Singularities 
--- Symposium I,} 
Lecture Notes in Mathematics {\bf 192},
Springer, Berlin, 1971, 165--168.

\bibitem{To06} V.A. Toponogov,
{\it Differential geometry of curves and surfaces,}
Birkh\"{a}user, Boston, 2006.
ISBN: 978-0-8176-4384-3; 0-8176-4384-2.
(Русский перевод: В.А. Топоногов,
{\it Дифференциальная геометрия кривых и поверхностей},
М.: Физматкнига, 2012. ISBN: 978-5-89155-213-5.)

\bibitem{Ho51}  H. Hopf,
{\it \"{U}ber Fl\"{a}chen mit einer Relation zwischen den 
Hauptkr\"{u}mmungen},
Math. Nachr. {\bf 4} (1951), 232--249.

\bibitem{Ho83}  H. Hopf,
{\it Differential geometry in the large},
Lecture Notes in Mathematics {\bf 1000}, 
Springer, Berlin, 1983.
      
\bibitem{Ya82}  S.-T. Yau,
{\it Problem section of the seminar in differential geometry at Tokyo,} 
Semin. differential geometry, Ann. Math. Stud. {\bf 102} (1982),
669--706.

\end{thebibliography}
\end{document}